\documentclass[12pt]{article}
\newtheorem{thm}{Theorem}

\newtheorem{cor}[thm]{Corollary}
\newtheorem{conj}[thm]{Conjecture}

\newcommand{\QED}{\hfill$\fbox{}$\vspace*{3mm}}
\newcommand{\Aut}{\mbox{Aut}}
\newcommand{\subgp}[1]{\langle{#1}\rangle}
\newcommand{\beeq}{\begin{eqnarray*}}
\newcommand{\eneq}{\end{eqnarray*}}
\newcommand{\proof}{\noindent{\it Proof.\hspace{4mm}}}
\newcommand{\example}{\noindent{\bf Examples.\hspace{4mm}}}
\newcommand{\qfd}{\hfill $\fbox{}$\vspace{4mm}}\def\newpic#1{%
\def\emline##1##2##3##4##5##6{%
\put(##1,##2){\special{em:point #1##3}}%
\put(##4,##5){\special{em:point #1##6}}%
\special{em:line #1##3,#1##6}}}
\newpic{}
\def\emline#1#2#3#4#5#6{%
\put(#1,#2){\special{em:moveto}}%
\put(#4,#5){\special{em:lineto}}}
\def\newpic#1{}
\newcommand{\Z}{\hbox{\bf Z}}
\newcommand\RR{\hbox{I\kern-.2em\hbox{R}}}
\newcommand\sRR{{\sl \hbox{I\kern-.2em\hbox{R}}}}\newcommand{\x}{\underline}
\newcommand{\y}{\underline}
\title{Quasiperfect domination in triangular lattices}
\author{Italo J. Dejter
\\ University of Puerto Rico \\ Rio Piedras, PR 00931-3355 \\ ijdejter@uprrp.edu
}\date{}

\begin{document}
\maketitle

\begin{abstract} A vertex subset $S$ of a graph $G$ is a perfect
(resp. qua\-si\-per\-fect) dominating set in $G$ if each vertex $v$
of $G\setminus S$ is adjacent to only one vertex ($d_v\in\{1,2\}$
vertices) of $S$. Perfect and qua\-si\-per\-fect
dominating sets in the regular tessellation graph of Schl\"afli
symbol $\{3,6\}$ and in its toroidal quotients are investigated, yielding
the classification of their perfect dominating sets and most of
their quasiperfect dominating sets $S$ with induced components of
the form $K_{\nu}$, where $\nu\in\{1,2,3\}$ depends only on $S$.
\end{abstract}

\section{Introduction}

A vertex subset $S$ of a graph $G$ is a {\it perfect dominating set}, or PDS, in $G$ if each vertex $v$ of the complementary graph $G\setminus S$ of $S$ in $G$ is adjacent to exactly one vertex of $S$. If $G$ is a regular graph, then $G\setminus S$ is a regular graph, \cite{Dej,DD2}. If in addition $S$ is isolated in $G$, then $S$ is
said to be a {\it 1-perfect code}, \cite{K}, or {\it efficient
dominating set} in $G$, \cite{BBS}.

In the present work, the graph that spans the regular tessellation
of Schl\"afli symbol $\{3,6\}$ in an Euclidean plane $\Pi\equiv\RR^2$ (with six
equilateral triangles at each vertex, as in the center of Fig. 4/4 in page
24 of \cite{Fej}) is referred to as the {\it
triangular lattice} $\Delta$. We will visualize $\Pi$ as the set of points $(x_1,x_2,x_3)\in\RR^3$ such that $x_1+x_2+x_3=0$; the vertices of $\Delta$ will be indicated by these
$\Delta$-{\it coordinates} in $\RR^3$, as shown in Figure 1. Even though
no two $\Delta$-coordinates are orthogonal, they are
more manageable for our purposes than orthogonal
coordinates are. Notice, for example, that the positive coordinate
directions of $x_1$ and $x_2$ are separated by an angle of $60^\circ$.

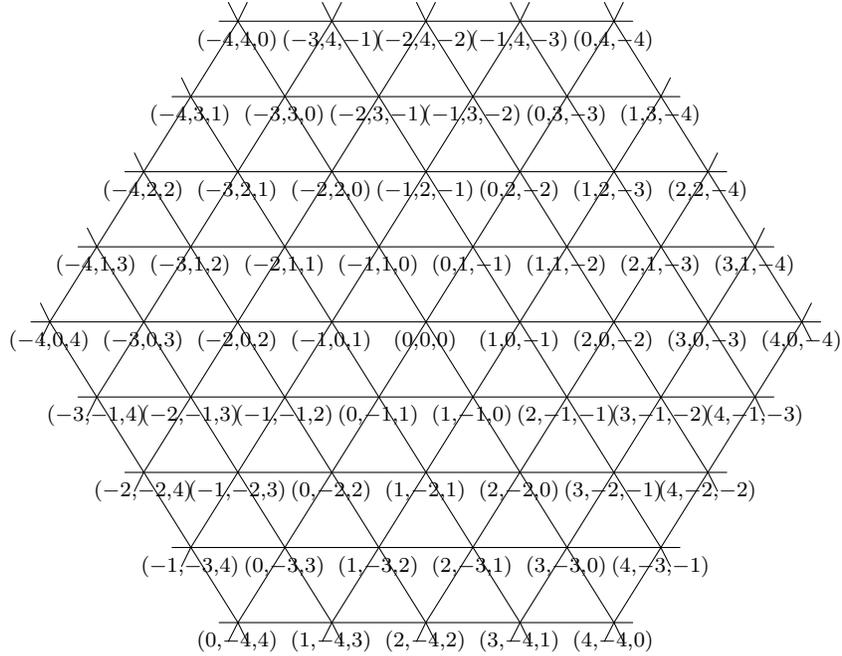
\begin{figure}
\unitlength=1.25mm
\special{em:linewidth 0.4pt}
\linethickness{0.4pt}
\begin{picture}(96.00,68.00)
\put(54.00,32.00){\makebox(0,0)[cc]{$_{(0,0,0)}$}}
\put(64.00,32.00){\makebox(0,0)[cc]{$_{(1,0,-1)}$}}
\put(74.00,32.00){\makebox(0,0)[cc]{$_{(2,0,-2)}$}}
\put(84.00,32.00){\makebox(0,0)[cc]{$_{(3,0,-3)}$}}
\put(94.00,32.00){\makebox(0,0)[cc]{$_{(4,0,-4)}$}}
\put(14.00,32.00){\makebox(0,0)[cc]{$_{(-4,0,4)}$}}
\put(24.00,32.00){\makebox(0,0)[cc]{$_{(-3,0,3)}$}}
\put(34.00,32.00){\makebox(0,0)[cc]{$_{(-2,0,2)}$}}
\put(44.00,32.00){\makebox(0,0)[cc]{$_{(-1,0,1)}$}}
\put(59.00,40.00){\makebox(0,0)[cc]{$_{(0,1,-1)}$}}
\put(69.00,40.00){\makebox(0,0)[cc]{$_{(1,1,-2)}$}}
\put(79.00,40.00){\makebox(0,0)[cc]{$_{(2,1,-3)}$}}
\put(89.00,40.00){\makebox(0,0)[cc]{$_{(3,1,-4)}$}}
\put(19.00,40.00){\makebox(0,0)[cc]{$_{(-4,1,3)}$}}
\put(29.00,40.00){\makebox(0,0)[cc]{$_{(-3,1,2)}$}}
\put(39.00,40.00){\makebox(0,0)[cc]{$_{(-2,1,1)}$}}
\put(49.00,40.00){\makebox(0,0)[cc]{$_{(-1,1,0)}$}}
\put(64.00,48.00){\makebox(0,0)[cc]{$_{(0,2,-2)}$}}
\put(74.00,48.00){\makebox(0,0)[cc]{$_{(1,2,-3)}$}}
\put(84.00,48.00){\makebox(0,0)[cc]{$_{(2,2,-4)}$}}
\put(24.00,48.00){\makebox(0,0)[cc]{$_{(-4,2,2)}$}}
\put(34.00,48.00){\makebox(0,0)[cc]{$_{(-3,2,1)}$}}
\put(44.00,48.00){\makebox(0,0)[cc]{$_{(-2,2,0)}$}}
\put(54.00,48.00){\makebox(0,0)[cc]{$_{(-1,2,-1)}$}}
\put(69.00,56.00){\makebox(0,0)[cc]{$_{(0,3,-3)}$}}
\put(79.00,56.00){\makebox(0,0)[cc]{$_{(1,3,-4)}$}}
\put(29.00,56.00){\makebox(0,0)[cc]{$_{(-4,3,1)}$}}
\put(39.00,56.00){\makebox(0,0)[cc]{$_{(-3,3,0)}$}}
\put(49.00,56.00){\makebox(0,0)[cc]{$_{(-2,3,-1)}$}}
\put(59.00,56.00){\makebox(0,0)[cc]{$_{(-1,3,-2)}$}}
\put(74.00,64.00){\makebox(0,0)[cc]{$_{(0,4,-4)}$}}
\put(34.00,64.00){\makebox(0,0)[cc]{$_{(-4,4,0)}$}}
\put(44.00,64.00){\makebox(0,0)[cc]{$_{(-3,4,-1)}$}}
\put(54.00,64.00){\makebox(0,0)[cc]{$_{(-2,4,-2)}$}}
\put(64.00,64.00){\makebox(0,0)[cc]{$_{(-1,4,-3)}$}}
\put(74.00,0.00){\makebox(0,0)[cc]{$_{(4,-4,0)}$}}
\put(34.00,0.00){\makebox(0,0)[cc]{$_{(0,-4,4)}$}}
\put(44.00,0.00){\makebox(0,0)[cc]{$_{(1,-4,3)}$}}
\put(54.00,0.00){\makebox(0,0)[cc]{$_{(2,-4,2)}$}}
\put(64.00,0.00){\makebox(0,0)[cc]{$_{(3,-4,1)}$}}
\put(69.00,8.00){\makebox(0,0)[cc]{$_{(3,-3,0)}$}}
\put(79.00,8.00){\makebox(0,0)[cc]{$_{(4,-3,-1)}$}}
\put(29.00,8.00){\makebox(0,0)[cc]{$_{(-1,-3,4)}$}}
\put(39.00,8.00){\makebox(0,0)[cc]{$_{(0,-3,3)}$}}
\put(49.00,8.00){\makebox(0,0)[cc]{$_{(1,-3,2)}$}}
\put(59.00,8.00){\makebox(0,0)[cc]{$_{(2,-3,1)}$}}
\put(64.00,16.00){\makebox(0,0)[cc]{$_{(2,-2,0)}$}}
\put(74.00,16.00){\makebox(0,0)[cc]{$_{(3,-2,-1)}$}}
\put(84.00,16.00){\makebox(0,0)[cc]{$_{(4,-2,-2)}$}}
\put(24.00,16.00){\makebox(0,0)[cc]{$_{(-2,-2,4)}$}}
\put(34.00,16.00){\makebox(0,0)[cc]{$_{(-1,-2,3)}$}}
\put(44.00,16.00){\makebox(0,0)[cc]{$_{(0,-2,2)}$}}
\put(54.00,16.00){\makebox(0,0)[cc]{$_{(1,-2,1)}$}}
\put(59.00,24.00){\makebox(0,0)[cc]{$_{(1,-1,0)}$}}
\put(69.00,24.00){\makebox(0,0)[cc]{$_{(2,-1,-1)}$}}
\put(79.00,24.00){\makebox(0,0)[cc]{$_{(3,-1,-2)}$}}
\put(89.00,24.00){\makebox(0,0)[cc]{$_{(4,-1,-3)}$}}
\put(19.00,24.00){\makebox(0,0)[cc]{$_{(-3,-1,4)}$}}
\put(29.00,24.00){\makebox(0,0)[cc]{$_{(-2,-1,3)}$}}
\put(39.00,24.00){\makebox(0,0)[cc]{$_{(-1,-1,2)}$}}
\put(49.00,24.00){\makebox(0,0)[cc]{$_{(0,-1,1)}$}}
\emline{34.00}{2.00}{1}{74.00}{2.00}{2}
\emline{74.00}{2.00}{3}{94.00}{34.00}{4}
\emline{94.00}{34.00}{5}{74.00}{66.00}{6}
\emline{74.00}{66.00}{7}{34.00}{66.00}{8}
\emline{34.00}{66.00}{9}{14.00}{34.00}{10}
\emline{14.00}{34.00}{11}{34.00}{2.00}{12}
\emline{29.00}{10.00}{13}{79.00}{10.00}{14}
\emline{79.00}{10.00}{15}{44.00}{66.00}{16}
\emline{44.00}{66.00}{17}{19.00}{26.00}{18}
\emline{19.00}{26.00}{19}{89.00}{26.00}{20}
\emline{89.00}{26.00}{21}{64.00}{66.00}{22}
\emline{64.00}{66.00}{23}{29.00}{10.00}{24}
\emline{24.00}{18.00}{25}{84.00}{18.00}{26}
\emline{84.00}{18.00}{27}{54.00}{66.00}{28}
\emline{54.00}{66.00}{29}{24.00}{18.00}{30}
\emline{74.00}{66.00}{31}{34.00}{2.00}{32}
\emline{14.00}{34.00}{33}{94.00}{34.00}{34}
\emline{34.00}{66.00}{35}{74.00}{2.00}{36}
\emline{29.00}{58.00}{37}{64.00}{2.00}{38}
\emline{64.00}{2.00}{39}{89.00}{42.00}{40}
\emline{89.00}{42.00}{41}{19.00}{42.00}{42}
\emline{19.00}{42.00}{43}{44.00}{2.00}{44}
\emline{44.00}{2.00}{45}{79.00}{58.00}{46}
\emline{79.00}{58.00}{47}{29.00}{58.00}{48}
\emline{24.00}{50.00}{49}{84.00}{50.00}{50}
\emline{84.00}{50.00}{51}{54.00}{2.00}{52}
\emline{54.00}{2.00}{53}{24.00}{50.00}{54}
\emline{14.00}{34.00}{55}{13.00}{32.00}{56}
\emline{14.00}{34.00}{57}{13.00}{36.00}{58}
\emline{14.00}{34.00}{59}{12.00}{34.00}{60}
\emline{19.00}{26.00}{61}{18.00}{24.00}{62}
\emline{19.00}{26.00}{63}{17.00}{26.00}{64}
\emline{24.00}{18.00}{65}{23.00}{16.00}{66}
\emline{24.00}{18.00}{67}{22.00}{18.00}{68}
\emline{29.00}{10.00}{69}{28.00}{8.00}{70}
\emline{29.00}{10.00}{71}{27.00}{10.00}{72}
\emline{34.00}{2.00}{73}{33.00}{0.00}{74}
\emline{34.00}{2.00}{75}{32.00}{2.00}{76}
\emline{35.00}{0.00}{77}{34.00}{2.00}{78}
\emline{44.00}{2.00}{79}{43.00}{0.00}{80}
\emline{45.00}{0.00}{81}{44.00}{2.00}{82}
\emline{54.00}{2.00}{83}{53.00}{0.00}{84}
\emline{55.00}{0.00}{85}{54.00}{2.00}{86}
\emline{64.00}{2.00}{87}{63.00}{0.00}{88}
\emline{65.00}{0.00}{89}{64.00}{2.00}{90}
\emline{74.00}{2.00}{91}{73.00}{0.00}{92}
\emline{75.00}{0.00}{93}{74.00}{2.00}{94}
\emline{76.00}{2.00}{95}{74.00}{2.00}{96}
\emline{19.00}{42.00}{97}{18.00}{44.00}{98}
\emline{19.00}{42.00}{99}{17.00}{42.00}{100}
\emline{24.00}{50.00}{101}{23.00}{52.00}{102}
\emline{24.00}{50.00}{103}{22.00}{50.00}{104}
\emline{29.00}{58.00}{105}{28.00}{60.00}{106}
\emline{29.00}{58.00}{107}{27.00}{58.00}{108}
\emline{34.00}{66.00}{109}{32.00}{66.00}{110}
\emline{76.00}{66.00}{111}{74.00}{66.00}{112}
\emline{80.00}{8.00}{113}{79.00}{10.00}{114}
\emline{81.00}{10.00}{115}{79.00}{10.00}{116}
\emline{85.00}{16.00}{117}{84.00}{18.00}{118}
\emline{86.00}{18.00}{119}{84.00}{18.00}{120}
\emline{90.00}{24.00}{121}{89.00}{26.00}{122}
\emline{91.00}{26.00}{123}{89.00}{26.00}{124}
\emline{95.00}{32.00}{125}{94.00}{34.00}{126}
\emline{96.00}{34.00}{127}{94.00}{34.00}{128}
\emline{95.00}{36.00}{129}{94.00}{34.00}{130}
\emline{91.00}{42.00}{131}{89.00}{42.00}{132}
\emline{90.00}{44.00}{133}{89.00}{42.00}{134}
\emline{86.00}{50.00}{135}{84.00}{50.00}{136}
\emline{85.00}{52.00}{137}{84.00}{50.00}{138}
\emline{81.00}{58.00}{139}{79.00}{58.00}{140}
\emline{80.00}{60.00}{141}{79.00}{58.00}{142}
\emline{75.00}{68.00}{143}{74.00}{66.00}{144}
\emline{74.00}{66.00}{145}{73.00}{68.00}{146}
\emline{65.00}{68.00}{147}{64.00}{66.00}{148}
\emline{64.00}{66.00}{149}{63.00}{68.00}{150}
\emline{55.00}{68.00}{151}{54.00}{66.00}{152}
\emline{54.00}{66.00}{153}{53.00}{68.00}{154}
\emline{45.00}{68.00}{155}{44.00}{66.00}{156}
\emline{44.00}{66.00}{157}{43.00}{68.00}{158}
\emline{35.00}{68.00}{159}{34.00}{66.00}{160}
\emline{34.00}{66.00}{161}{33.00}{68.00}{162}
\end{picture}
\caption{Simplex coordinates in $\Delta$}
\end{figure}

Suggested by A. Delgado, the problem of characterizing the PDSs in $\Delta$ and its toroidal quotients is completed in Sections 2-3: Up to symmetry, there is just one 1-perfect code in $\Delta$. In order to broaden the subject, the PDS condition is relaxed by defining a {\it quasiperfect dominating set}, or QPDS, in a graph $G$ as a vertex subset $S$ of $G$ for which each vertex $v$ of $G\setminus S$ is adjacent to $d_v\in\{1,2\}$ vertices of $S$.

Section 4 deals with constant $d_v=2$; the rest of the paper is devoted to characterizing, in $\Delta$ and its toroidal quotients, the QPDSs $S$ with induced components of the form $K_\nu$, for fixed $\nu=\nu(S)\in\{1,2,3\}$, (only dependent on $S$). This is
left open in $\Delta$ for a few unknown possibilities, in Theorem 13 and Conjecture 15. Theorems 1, 6 and 22 characterize all the remaining cases; they have such $K_\nu$'s with a minimum graph distance $\delta=3$, in contrast with one in Section 4, for which $\delta=2$, (see left of Figure 3).

An ordered triple of pairwise adjacent vertices of $\Delta$, as for example $T_0=((0,0,0),(1,0,-1),(0,1,-1))$, will be said to be an {\it ordered triangle}. There is a bijection $\Theta$ from the group ${\mathcal A}(\Delta)$ of automorphisms of $\Delta$ onto the
collection of ordered triangles. In fact, each $\phi\in{\mathcal A}(\Delta)$
is determined by assigning $T_0$ to an ordered triangle $\Theta(\phi)=\phi(T_0)$ with the vertices of $T_0$ sent bijectively (in one of six fashions) onto the vertices of $\Theta(\phi)=\phi(T_0)$.
This yields that ${\mathcal A}(\Delta)$ is a semidirect product of the cyclic group
$D_{12}$ (of symmetries of $\Delta$ around $(0,0,0)$) and the group
$\Z^2$ (of parallel translations of $\Delta$ expressed, say, in the first two
$\Delta$-coordinates: $x_1,x_2$). This semidirect product,
$D_{12}\times_\phi\Z^2$, is given in the Cartesian set product
$D_{12}\times\Z^2$ via the homomorphism
$\phi:\Z^2\rightarrow{\mathcal A}(D_{12})$ defined by
$\phi(m)(r)=mrm^{-1}=mr(-m)$, with multiplication set by
$(r,m)(s,n)=(r\phi(m)(s),mn)$, where $r,s\in D_{12}$ and
$(m,n)\in\Z^2$.

\section{Perfect domination in $\Delta$}

\begin{figure}
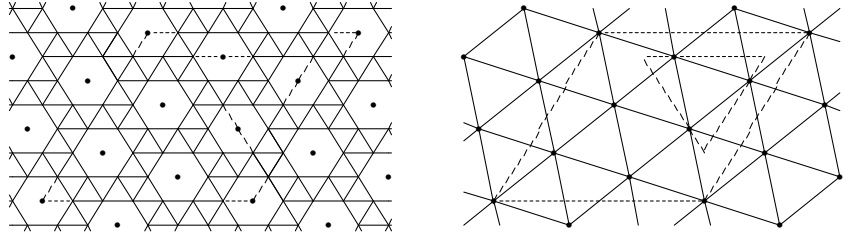

\unitlength=0.40mm
\special{em:linewidth 0.4pt}
\linethickness{0.4pt}

\caption{Complement of 1-Perfect code $S$ in $\Delta$ and its graph
$M(S)$}
\end{figure}

The {\it minimum distance graph} $M(S)$ of a 1-perfect code $S$ in $\Delta$ is
the graph whose vertex set is $S$ and whose edges are traced between
those pairs of vertices that realize the minimum graph distance, 3, of $S$ in
$\Delta$. On the left of
Figure 2, such an $S$ and its 5-regular
complement $\Delta\setminus S$ in $\Pi$ are represented, (coinciding with the $(3,3,3,3,6)$-tessellation of Fig. 5/1 in page
42 of \cite{Fej}); its $M(S)$ is shown on the right of Figure
2 (superposable, vertex by corresponding vertex, over the left of Figure 2).

It is convenient to
represent the line graph $L(S)$ of $M(S)$ as a plane graph with its vertices as the midpoints
of the edges of $M(S)$. By properly coloring the faces of $L(S)$
with seven colors, it can be seen that an upper bound on the
chromatic number $\chi$ of the graph of $\Pi$ with a forbidden fixed
Euclidean distance is 7, (which can be realized via the chromatic number, 7, of $M(S)$;  actually, Thomassen \cite{Tho} showed that $\chi=7$).

The complement of $\Delta\setminus S$ in $\Pi$ is the disjoint
union of 2-dimensional connected components whose closures are
regular hexagons and equilateral triangles, both of unit side
lengths. The Euclidean distance between such hexagons (resp.
triangles) has a lower bound of $\frac{\sqrt{3}}{2}$ (0). The
collection of these hexagons and triangles will be denoted by
${\mathcal P}(S)$. An $n$-{\it hole} of a graph $G$ is
an induced $n$-cycle of $G$ (no chords, i.e. extra edges). Then
the boundary of each hexagon $h$ in ${\mathcal P}(S)$ forms a
6-{\it hole} of $\Delta\setminus S$. Alternatively, we can say that
$h$ is {\it bordered} by a 6-hole.

A superposed plane representation of both $L(S)$ and
$\Delta\setminus S$ can be found on the lower part of Fig. 5/2 in
page 46 of \cite{Fej}, where the vertices of $L(S)$ are the
centroids of those triangles in ${\mathcal P}(S)$ whose vertices are
shared by three different hexagons of ${\mathcal P}(S)$. Each edge
of $L(G)$ is of length $\sqrt{\frac{7}{3}}$ and perpendicular to the
edge of $M(G)$ that crosses it. The edge perpendicular to each edge
$e$ of $M(G)$ is realized by the middle third of the segment joining
the pair of vertices opposite to $e$ in the triangles incident to
$e$.

\begin{thm}
Up to symmetry, there is exactly one proper {\rm PDS} $S$ in $\Delta$.
Moreover, $S$ is a 1-perfect code and its minimum graph distance is $\delta=3$.
\end{thm}

\proof We assert that a proper PDS $S$ in $\Delta$ is isolated. For
otherwise there would be two adjacent vertices $v$ and $w$ in $S$
having a common neighbor in $\Delta\setminus S$, a contradiction. In
fact, there exists exactly two proper PDS $S$ in $\Delta$, up to
parallel translations and rotations of $\Pi$; their complementary graphs
$\Delta\setminus S$ are {\it enantiomorphic} (that is:
mirror images of each other, say by reflection on the $x_1$-axis), so they are isomorphic. One of these PDSs is the one depicted on the left of Figure 2. \qfd

Let $S$ be as in Theorem 1. An ordered pair of adjacent vertices of
$\Delta$ is said to be an {\it arc}. Let $S$ be a PDS in $\Delta$.
Then there exists a bijection $\theta$ from
${\mathcal A}(\Delta\setminus S)$ onto the collection of arcs of 6-holes of
$\Delta\setminus S$. Assume that $\Delta\setminus S$ has a 6-hole containing the arc $A_0=((0,-1,1),(1,-1,0))$. Then
any $\psi\in{\mathcal A}(\Delta\setminus S)$
is determined by assigning $A_0$ to $\theta(\psi)=\psi(A_0)$. This allows us to see ${\mathcal A}(\Delta\setminus S)$ as a semidirect product of $\Z_6$ and $\Z^2$, where the generators
$(1,0)$ and $(0,1)$ of $\Z^2$ are sent onto the parallel translations of $\Delta\setminus S$ along respective vectors $(3,-1-2)$ and $(1,2,-3)$.
However, ${\mathcal A}(M(S))={\mathcal A}(L(S))$ is isomorphic to ${\mathcal A}(\Delta)$.

Each vertex $x$ of $S$ excludes its six neighbors in $\Delta$ from
$S$, which determine a regular 6-hole $S(x)$ in
$\Delta\setminus S$; each edge of $S(x)$ is adjacent to a
triangle in $\Delta\setminus S$, so there are six triangles in
$\Delta\setminus S$ adjacent to $S(x)$ by means of an edge each.
If a triangle $t$ of $\Delta\setminus S$ is not adjacent to any $S(x)$, then $t$ is incident to three 6-holes $S(x),S(y),S(z)$, where $x,y,z\in S$ are at unit distance from $t$.
Any such $t$ will be called a $\Z_3$-triangle;
(rotations around the barycenter of $t$ offer three automorphisms of
$\Delta\setminus S$ composing a subgroup $\Z_3$ of
${\mathcal A}(\Delta\setminus S)$).
So, each vertex of a $\Z_3$-triangle is incident to just one 6-hole
in $\Delta\setminus S$.

\begin{cor} $V(\Delta)$ admits a partition into seven translated
copies of $S$ at unit distance from each other. Moreover,
the total number of {\rm PDS}s isomorphic to $S$ in $\Delta$ is 14.
These 14 {\rm PDS}s compose two enantiomorphic partitions of $\Delta$ into 1-perfect codes. \end{cor}

\proof Given a fixed $S(x)$ in $\Delta\setminus S$, the parallel translations of $S$ taking $x$ to its neighbors in $S(x)$ give the remaining six members of the partition in the first assertion of the corollary. This yields seven of the 14 copies of $S$ mentioned in the second assertion. The other seven copies of $S$ are obtained similarly from a reflected copy of $S$ in $\Delta$
through the reflection, say, on the $x_1$-axis. \qfd

Note that the inner dotted lines on the left of Figure 2 form
an example of the largest 1-perfect code existing on a triangular
subgraph of $\Delta$.

\section{PDSs in toroidal triangular lattices}

Let us express the vertices of $\Delta$ by means of their first two $\Delta$-coordinates: $x_1,x_2$. The subgraph $\Delta'$ of $\Delta$ spanned by the edges with constant $x_1$ or $x_2$  has a quotient Cartesian product $C_m\times C_n$ of cycles of lengths $m,n\in\Z$ with $3\leq m$ and $3\leq n$. From such $C_m\times C_n$, which is a toroidal graph (i.e. embeddable into the flat torus ${\mathcal T}$), we obtain a {\it triangular lattice} graph $\Delta_{
m,n}$ by adjoining to it the {\it anti-diagonal edges} $\{(i,j+1),(i+1,j)\}$ of the {\it elementary $4$-cycles} $((i,j),(i+1,j),(i+1,j+1),(i,j+1))$ of $C_m\times C_n$, where $0\leq i<m$, $0\leq j<n$ and additions are taken mod $m,n$,  respectively.

There are 1-perfect codes $S$ in the graphs $\Delta_{7k,7\ell}$, for
$0<k,\ell\in\Z$, with cardinality
$\frac{1}{7}|V(\Delta_{7k,7\ell})|=7k\ell$. An example of such an $S$
in $\Delta_{7,7}$ can be visualized from either side of Figure 2,
where the outer dotted lines delineate the boundary of a rhomboidal
pattern $\Upsilon$, which is a cutout of ${\mathcal T}$.

\begin{cor}
The vertex set of a graph $\Delta_{m,n}$ has a 1-perfect code
partition ${\mathcal S}=\{S^0=S,S^1,\ldots,S^6\}$ if and only if $m$
and $n$ are multiples of 7. Each component code $S^i$
is a translate of $S$ and has cardinality $mn/7$.
\end{cor}

\proof Using that $\Upsilon$
is a minimal cutout for a toroidal embedding of
$\Delta$, it can be seen that the natural projection
$\Delta\rightarrow\Delta_{m,n}$ maps a partition as in Corollary 2
onto a partition as claimed.
By allowing larger rhomboidal cutouts formed by $m/7$ horizontally
contiguous copies of $\Upsilon$ times $n/7$ diagonally contiguous copies
of $\Upsilon$, the existence of the claimed partitions is
ensured. Other values of $m$ and $n$ are clearly incompatible with
the existence of isolated PDSs in corresponding graphs
$\Delta_{m,n}$. \qfd

The group ${\mathcal A}(\Delta_{m,n})$ can be expressed with $\Z^2$ in the remark of the
last paragraph in Section 1 replaced by $\Z_m\times\Z_n$. It follows that
${\mathcal A}(\Delta_{m,n}\setminus S)$ is a semidirect product of the
rotation group $\Z_6$ and the translation group
$\Z_m\times\Z_{n/7}$, or also of $\Z_6$ and $\Z_{m/7}\times\Z_n$.

\section{Semiperfect domination in $\Delta$ and $\Delta_{m,n}$}

Let $G$ be a graph and let $S$ be a QPDS in $G$. If each vertex of
$G\setminus S$ is adjacent to {\it exactly two} vertices of $S$, then $S$ is a {\it semiperfect dominating set}, or SPDS, in $G$.

\begin{thm}
Up to symmetry, there are exactly two {\rm SPDS}s $S$ in $\Delta$.  One of them has  minimum graph distance $\delta$ between induced components of $S$  equal to 2, while the other has $\Delta$ equal to 3.
\end{thm}

\proof An isolated (resp. non-isolated) SPDS $S$ in $\Delta$ with $\delta=2$ (3) and induced components spanning
affine spaces of dimension 0 (1) in $\Pi$ is shown on the left (center) of Figure 3, with
connected (disconnected) 4-regular complement $\Delta\setminus S$ in $\Pi$; the two edges at each $v\in V(\Delta\setminus S)$ are
separated by an angle of $180^\circ$ ($60^\circ$) and their endvertices in $S$ are in different (the same) component(s) of $\Delta[S]$. Up to symmetry, no other SPDSs exist in $\Delta$.\qfd

The graph $\Delta\setminus S$ on the  left of Figure 3 coincides with the $(3,6,3,6)$-tessellation
of Fig. 5.1 in page 42 of \cite{Fej}. Its
automorphism group, $D_{12}\times_\phi(2\Z)^2$ (isomorphic to
$D_{12}\times_\phi\Z^2$), is embedded naturally into
${\mathcal A}(\Delta)=D_{12}\times_\phi\Z^2$.

\begin{figure}
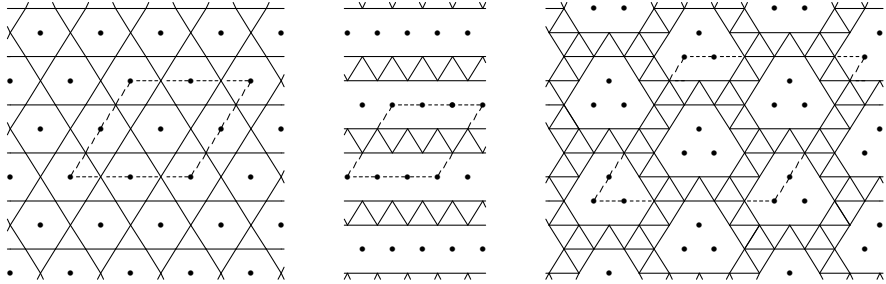

\unitlength=0.40mm
\special{em:linewidth 0.4pt}
\linethickness{0.4pt}

\caption{Complements of QPDSs in $\Delta$ in Sections 4 and 5}
\end{figure}

Now, let $S$ be as in the center of Figure 3. The automorphism group
${\mathcal A}(\Delta\setminus S)$ is the semidirect product of the doubly
reflective group $\Z_2^2$ (whose reflection axes are the $x_1$-axis and the line that passes through $(0,0,0)$ and $(-1,2,-1)$) and
the translation group $\Z^2$, with its generators
$(1,0)$ and $(0,1)$ of $\Z^2$  sent onto the parallel translations of $\Delta\setminus S$ along respective vectors $(1,0,-1)$ and $(0,3,-3)$.

\begin{cor}{\bf(A)} For $1\leq m,n\in\Z$, $\Delta_{2m,2n}$ contains
an isolated {\rm SPDS} $S$, with connected induced subgraph
$\Delta_{m,n}[S]$. {\bf(B)} For $3\leq m,n\in\Z$, there exists a
graph $\Delta_{m,n}$ containing a non-isolated SPDS $S$ if and only
if either $3|m$ or $3|n$. Accordingly, the number of induced
components of $S$ in $\Delta_{m,n}$, or components of
$\Delta_{m,n}[S]$, is either $\frac{m}{3}$ or $\frac{n}{3}$, so each
such component is either an $n$-cycle or an $m$-cycle.
\end{cor}

\proof Examples of SPDSs as in items (A) and (B) can be
visualized from the left and center of Figure 3, respectively, where
the dotted lines delineate the boundaries of rhomboidal cutouts of ${\mathcal T}$.
In the first case of (B), for example,
there are SPDSs in the graphs $\Delta_{k,3\ell}$ for $3\leq k\in\Z$ and
$1\leq\ell\in\Z$. These SPDSs have cardinality
$\frac{1}{3}|V(\Delta_{k,3\ell})|=3k\ell$. \qfd

\section{$K_3$-quasiperfect domination in $\Delta$ and $\Delta_{m,n}$}

Let $H$ be a subgraph of a graph $G$. A QPDS $S$ of $G$ is
an $H$-{\it quasiperfect dominating set}, or $H$-QPDS, in $G$ if all
the induced components of $S$ in $G$, that is all the components of
$G[S]$, are isomorphic to $H$.

\begin{thm} Up to symmetry, there is exactly one $K_3$-{\rm QPDS} in $\Delta$.  The minimum graph distance between $K_3$'s induced in $S$ is $\delta=3$.
\end{thm}

\proof If $S$ is a $K_3$-QPDS in $\Delta$, then the complement of
$\Delta\setminus S$ in $\Pi$ is the disjoint union of
2-dimensional connected components whose closures are equilateral
triangles and equiangular semisymmetric hexagons (the elements of a set ${\mathcal P}(S)$, as in Section 2), illustrated on
the left of Figure 3. The Euclidean distance between such
hexagons (resp. triangles) has a lower bound of $\frac{\sqrt{3}}{2}$
(0). Notice that the boundary of each such hexagon
constitutes a 9-hole in $\Delta\setminus S$. Clearly, $\delta=3$.\qfd

The automorphism group ${\mathcal A}(\Delta\setminus S)$ for the $K_3$-QPDS
$S$ of $\Delta$ in the proof of Theorem 6 is a semidirect product of
the group $S_3$ of symmetries of a fixed equiangular semiregular
hexagon as above and the group $\Z^2$, where the generators
$(1,0)$ and $(0,1)$ of $\Z^2$ are sent onto the parallel translations of
$\Delta\setminus S$ along respective vectors $(4,-2,-2)$ and $(2,2,-4)$.

From Theorem 6 and the left of Figure 3, we see that for
$0<k,\ell\in\Z$ there are $K_3$-QPDSs in toroidal graphs
$\Delta_{6k,6\ell}$ and that they have cardinality
$\frac{1}{9}|V(\Delta_{6k,6\ell})|=3k\ell$. An example of such a
$K_3$-QPDS in a toroidal graph $\Delta_{6,6}$ can be visualized by
identifying adequately the dotted lines on the left of Figure 3.

\begin{cor}
There exists a toroidal graph $\Delta_{m,n}$ containing a $K_3$-{\rm QPDS}
$S$ if and only if $m$ and $n$ are multiples of 6. The number of
induced components of such an $S$ in $\Delta_{m,n}$, that is the
components of $\Delta_{m,n}[S]$, is $\frac{mn}{12}$.
\end{cor}

\proof The argument here resembles that of the proof of Corollary 3 but for just
one $K_3$-QPDS, because there is no partition in this case. \qfd

The group ${\mathcal A}(\Delta\setminus S)$ of the $K_3$-QPDS $S$ in
$\Delta_{m,n}$ of Theorem 6 is a semidirect product of
$\Z_m\times\Z_{n/3}$ and $S_3$, or also of $\Z_{m/3}\times\Z_n$ and
$S_3$.

\section{$K_2$-quasiperfect domination in $\Delta$ and $\Delta_{m,n}$}

Let $S$ be a $K_2$-QPDS in $\Delta$. The complement of
$\Delta\setminus S$ in $\Pi$ is the disjoint union of
2-dimensional connected components whose closures are equilateral
triangles and elongated hexagons, these containing each a unique
induced component $K_2$ of $S$ in $\Delta$ and thus a unique induced
edge. The Euclidean distances between such hexagons (resp.
triangles) has a lower bound of $\frac{\sqrt{3}}{2}$ (0).
The collection of these triangles and hexagons will be denoted by
${\mathcal P}(S)$. Notice that the boundary of each elongated
hexagon in ${\mathcal P}(S)$ here constitutes an 8-hole in
$\Delta\setminus S$. It is clear that the minimum graph distance between those components $K_2$ is $\delta=3$.

\begin{figure}
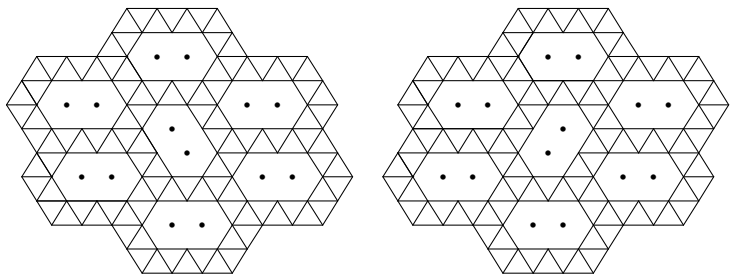

\unitlength=0.40mm
\special{em:linewidth 0.4pt}
\linethickness{0.4pt}

\caption{Examples of hexagon types in complements of $K_2$-QPDSs in $\Delta$}
\end{figure}

An hexagon in ${\mathcal P}(S)$ is said to be of {\it type}
$i\in I_3=\{1,2,3\}$ if the vertices of $S$ it contains share the
$i$-th $\Delta$-coordinate, $x_i$. Figure 4
contains on its left (resp. right) side an hexagon of type 2 (3) bordered by equilateral triangles separating it from six hexagons
of type 1. In view of this, we define an auxiliary graph
$\Gamma(S)$ whose vertices are the components of $\Delta[S]$, or of
some $\Delta_{m,n}[S]$; two such vertices are adjacent $\Gamma(S)$ if and only if the induced components they represent are at a graph distance of 3, the minimum attainable graph distance
separating the components of $\Delta[S]$ in $\Delta$. Alternatively, we may
consider the vertices of $\Gamma(S)$ as the hexagons of types $\in I_3$, any two of them adjacent if and only if the Euclidean
distance between them is $\frac{\sqrt{3}}{2}$.

For the two portions of a $\Delta\setminus S$
in Figure 4, the corresponding portions of their associated
graphs $\Gamma(S)$ look like as in the leftmost tables below:
$$\begin{array}{|ccccc|ccccccc|ccccccc|}\hline
&   &   &   & & &   &   &   &   &   & & &   &   &   &   &   & \\
& _1  & ^1 & _1  & & & ^1 & _1  & ^1 & _1  & ^2 & & & ^2 & _1  & ^1 & _1  & ^2 & \\
& _1  & ^2 & _1  & & & ^1 & _1  & ^2 & _1  & ^1 & & & ^1 & _1  & ^2 & _1  & ^1 & \\
&   & ^1 &   & & & ^2 &   & ^1 &   & ^1 & & & ^2 &   & ^1 &   & ^2 & \\
\hline
&   &   &   & & &   &   &   &   &   & & &   &   &   &   &   & \\
& _1  & ^1 & _1  & & & ^3 & _1  & ^1 & _1  & ^1 & & & ^3 & _1  & ^1 & _1  & ^3 & \\
& _1  & ^3 & _1  & & & ^1 & _1  & ^3 & _1  & ^1 & & & ^1 & _1  & ^3 & _1  & ^1 & \\
&   & ^1 &   & & & ^1 &   & ^1 &   & ^3 & & & ^3 &   & ^1 &   & ^3 & \\
\hline
\end{array}$$
where the type of each vertex represents it, and two vertices so
represented are adjacent if they are vertically, diagonally or
anti-diagonally contiguous, meaning they are on a line at
$90^\circ$, $60^\circ$ or $120^\circ$, res\-pectively, from
the $x_1$-axis. The tables at right and far right of
each such table represent two possible extended portions of
$\Delta\setminus S$ in each case. We conclude that there are two
extreme possibilities for extending these portions to complements of
$K_2$-QPDSs in $\Delta$. One of them has an infinitely
extended table with either a unique infinite anti-diagonal that descends or
ascends with period $(12)$, as in the bisequence (i.e. doubly
infinite sequence) $(\ldots 1212\ldots)$, or a unique infinite
diagonal that descends or ascends with period $(13)$, as in the
bisequence $(\ldots 1313\ldots)$, with all the hexagon types that
are away from this anti-diagonal or diagonal being equal to 1, as in
the middle tables of the arrangement above. The other extreme
possibility is to have alternate diagonals, (resp. anti-diagonals)
of periods $(12)$ ($(13)$) and $(1)$, as on the rightmost
tables above. Between these two pairs of extreme situations lie all
the possible complements of $K_2$-QPDSs in $\Delta$ whose induced
components pertain to just two hexagon types, one of which is 1.

\begin{figure}
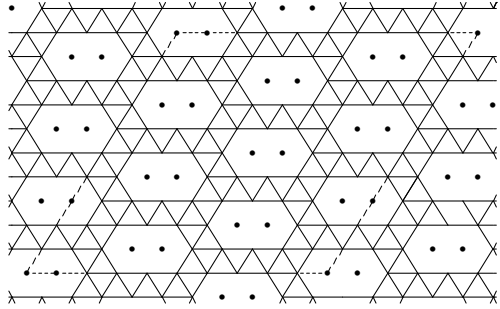

\unitlength=0.40mm
\special{em:linewidth 0.4pt}
\linethickness{0.4pt}

\caption{Complement of a type-parallel $K_2$-QPDS in $\Delta$}
\end{figure}

We consider now the $K_2$-QPDSs $S$ of $\Delta$ for which the hexagons in ${\mathcal P}(S)$ have constant type, for example type 1, as in Figure 5. The automorphism group
${\mathcal A}(\Delta\setminus S)$, for one such $S$, is the semidirect product of the doubly reflective group $\Z_2^2$ of a fixed hexagon of type 1 and the group $\Z^2$, with its
generators $(1,0)$ and $(0,1)$ sent onto the parallel translations of $\Delta\setminus S$ along respective vectors $(4,-2,-2)$ and $(3,1,-4)$.

Any $K_2$-QPDS in $\Delta$ for which the hexagons in ${\mathcal P}(S)$ have constant type
is said to be a {\it parallel} $K_2$-QPDS.

\begin{thm}
Up to symmetry, there is only one parallel $K_2$-{\rm QPDS} $S$ in $\Delta$.
\end{thm}

\proof Consider an hexagon of type 1 in ${\mathcal P}(S)$,
for an $S$ as in the statement, and delimited by one of its 8-holes.
By analyzing the possibilities of neighboring hexagons of type 1 in
the formation of such an $S$, it follows that the case in Figure 5 is the only possible one. \qfd

\begin{cor} For $0<k,\ell\in\Z$, there are parallel $K_2$-{\rm QPDS}s in the graphs
$\Delta_{10k,10\ell}$ with cardinality
$\frac{1}{10}|V(\Delta_{10k,10\ell})|=10k\ell$.\end{cor}

\begin{figure}
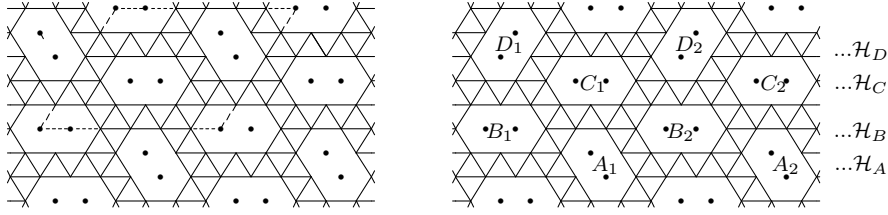

\unitlength=0.40mm
\special{em:linewidth 0.4pt}
\linethickness{0.4pt}

\caption{Complements of QPDSs in $\Delta$ for Theorems 11 and 16}
\end{figure}

\proof An example of such a parallel $K_2$-QPDS in a graph
$\Delta_{10,10}$ can be visualized by identifying adequately the
dotted lined in Figure 5.\qfd

\begin{cor}
The vertex set of a graph $\Delta_{m,n}$ has a parallel $K_2$-{\rm QPDS}
partition ${\mathcal S}=\{S^0=S,S^1,\ldots,S^9\}$ if and only if $m$
and $n$ are multiples of 10. Each component code $S^i$ of such a
partition ${\mathcal S}$ is a translate of $S$ and has
cardinality $mn/10$.
\end{cor}

\proof This follows from the cutout in Figure 5 and the ideas in the
proofs of Corollaries 3, 7 and 9.\qfd

Now, consider the case of $K_2$-QPDSs in $\Delta$, or some
$\Delta_{m,n}$, with just two types of hexagons, say types $i,j\in
I_3$, $i\neq j$, in which an hexagon of type $i$ is {\it surrounded}
by hexagons of type $j$: its adjacent hexagons in $\Gamma(S)$
are all of type $j$. In $\Delta$, this implies that each hexagon of
type $j$ is surrounded by three hexagons of type $i$ and three
alternated hexagons of type $j$. For example, the upper-right table
in the arrangement of six tables above can be extended so that each
hexagon of type 2 is surrounded by six hexagons of type 1.  We say
that a $K_2$-QPDS $S$ in $\Delta$ or some $\Delta_{m,n}$, with
exactly two hexagon types $i,j\in I_3$, is $(i,j)$-{\it surrounded} if
each of its hexagons of type $i$ is surrounded by hexagons of
type $j$. The following result is generalized in Theorem 16.

\begin{thm} Up to symmetry, there is only one
$(i,j)$-sur\-round\-ed $K_2$-{\rm QPDS} in $\Delta$, for each two types $i,j\in I_3$.
\end{thm}

\proof The $K_2$-QPDS in the left (resp. right) side of
Figure 4 extends to the sole existing $K_2$-QPDS as claimed, with $i=2$ (3) and $j=1$.
\qfd

\begin{cor} There exists a graph $\Delta_{m,n}$ containing a
$(2,1)$-sur\-round\-ed $K_2$-{\rm QPDS} $S$ if and only if $6|m$ and
$5|n$. The number of hexagons of type $2$, resp. $1$, for such an
$S$ in $\Delta_{m,n}$ is $\frac{mn}{30}$, resp. $\frac{2mn}{30}$.
\end{cor}

\newpage

\begin{figure}
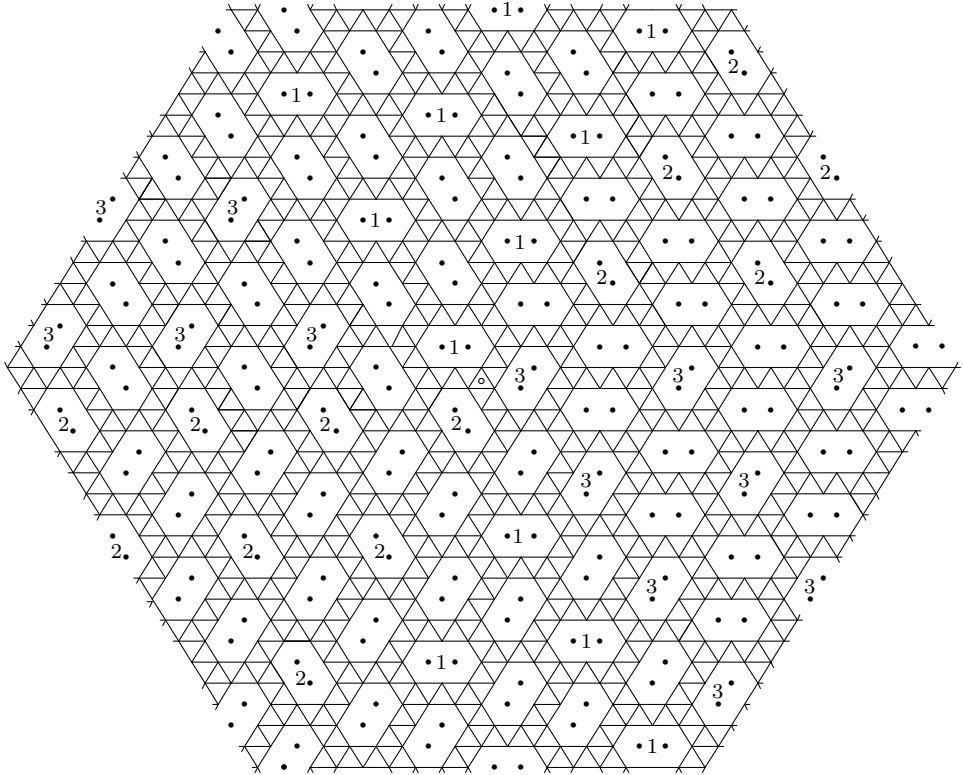

\unitlength=0.35mm
\special{em:linewidth 0.4pt}
\linethickness{0.4pt}

\caption{Complement of $K_2$-QPDS in $\Delta$ with $\Z_3$-symmetry}
\end{figure}

\proof This follows the ideas in the proofs of Corollaries 3, 7 and
9, with any of its smallest cases, for example the cutout on the
left representation of Figure 6, leading to a $K_2$-QPDS $S$ in
$\Delta_{6,5}$ that produces barely one hexagon of type 2, say $J$,
surrounded by only two hexagons of type 1, each adjacent to $J$
by means of three edges in $\Gamma(S)$.\qfd

\begin{thm} Up to symmetry, there are at least two $K_2$-{\it QPDS}s in $\Delta$ with three pairwise adjacent hexagons of different types, covering $I_3$.
\end{thm}

\proof Figure 7 can be extended to a $K_2$-QPDS $S$ for which only the hexagons of ${\mathcal P}(S)$ adjacent to more than three hexagons of a common type are labeled with its
type. The set of all the types for this figure has a disposition as
in the left frame of the following table, for which the mentioned
hexagons (adjacent to more than three hexagons of a common type) are
represented by underlined types, as a reference, and the hexagons
nearest to the center of $\Z_3$-symmetry are indicated each with a dot
over their type.
$$\begin{array}{|ccccccccccccc|ccccccccccccc|}\hline
&&&&&&&&&&&&&_1&&_3&&_{\x3}&&_1&&_3&&_1&&_{\y3}\\
&&&&_2&&_2&^1&_2&^1&_{\x2}&&&_3&^3&_1&^1&_3&^{\x3}&_3&^3&_1&^1&_{\y3}&^{\y3}&_{\y3}\\
&&&^2&_{\x1}&^2&_{\x1}&^2&_{\x1}&^1&_1&&&_1&^1&_3&^3&_1&^{\x3}&_{\x3}&^1&_{\y3}&^{\y3}&_{\y3}&^{\y3}&_{\y3}\\
&&_2&^2&_2&^2&_2&^2&_1&^{\x2}&_1&^{\x2}&&_1&^3&_1&^1&_3&^3&_{\x3}&^3&_2&^3&_2&^3&_2\\
&^{\x3}&_2&^{\x3}&_2&^{\x1}&_2&^{\x1}&_{\x2}&^1&_{\x2}&^1&&_{\y1}&^1&_3&^3&_1&^1&_3&^{\x3}&_2&^3&_2&^3&_2\\
_{\x3}&^2&_{\x3}&^2&_{\x3}&^2&_{\dot{\x1}}&^1&_1&^1&_1&^1&_1&_{\y1}&^{\y1}&_1&^1&_3&^3&_{\dot 1}&^{\x3}&_2&^3&_2&^3&_2\\
_{\x2}&^2&_{\x2}&^2&_{\x2}&^2&_{\dot{\x2}}&^{\dot{\x3}}&_1&^{\x3}&_1&^{\x3}&_1&_2&^{\y1}&_{\y1}&^3&_{\x1}&^{\x1}&_{\dot 2}&^{\dot 3}&_2&^3&_2&^3&_2\\
&^3&_3&^3&_3&^3&_3&^3&_{\x3}&^1&_{\x3}&^1&&_1&^2&_2&^1&_{\x1}&^1&_1&^{\x2}&_{\x2}&^{\x2}&_{\x2}&^{\x2}&_{\x2}\\
&^{\x2}&_3&^{\x2}&_3&^{\x2}&_3&^{\x1}&_3&^1&_1&^1&&_2&^1&_{\x1}&^{\x1}&_1&^2&_2&^1&_1&^2&_3&^2&_3\\
&&&^3&_{\x2}&^3&_{\x1}&^3&_{\x1}&^{\x3}&_1&^{\x3}&&_1&^2&_{\x1}&^1&_2&^1&_1&^2&_2&^{\y2}&_3&^2&_3\\
&&&^3&_3&^3&_3&^3&_3&^3&_{\x3}&&&_2&^{\x1}&_1&^1&_1&^2&_2&^1&_{\y2}&^{\y2}&_3&^2&_3\\
&&&&&&&^1&&^1&&&&&^{\x1}&&^2&&^1&&^2&&^{\y2}&&^2&\\\hline
\end{array}$$
The right frame of the table offers a (larger) corresponding
disposition of a second $K_2$-QPDS $S$ as in the statement, possessing
only six hexagons adjacent to at most three hexagons of a common
type. Among these six hexagons, there are three nearest to the
center of $\Z_3$-symmetry; each of them constitutes the initial
vertex of an infinite type-constant path $\lambda_i$ in
$\Gamma(S)$, ($i\in I_3$). Each $\lambda_i$,
whose vertices are indicated in the frame by underlined types, has
contiguous edges with two alternate directions. For example,
$\lambda_3$ has all its vertices labeled $\underline 3$, starts immediately
over the vertex labeled $\dot 3$ and proceeds first with a vertical edge, then
with a diagonal edge (in the upper-left direction), and
so on, alternatively. Note that the edge directions are
perpendicular to the edges of $\Delta[S]$ in the
hexagons that represent vertices of the path, in each $\lambda_i$.

The other three of the mentioned six hexagons are each the
initial vertex of two infinite type-constant paths sharing their
first edge: a rectilinear path $\eta_i$ (in the direction
perpendicular to the edge of $\Delta[S]$ in this initial hexagon)
and a zigzagging path $\xi_i$ ($i\in I_3$). The vertices in both
paths and between them have constant type, that we underline in the
frame (as we did for the vertices of $\lambda_i$, $i\in I_3)$. For
example toward the upper-right corner of the frame, a rectilinear
path $\eta_3$ of vertices labeled $\underline 3$ (showing five such
vertices in the frame) can be seen at an angle of $30^\circ$
over the $x_1$-axis, sharing its first edge with a
zigzagging path $\xi_3$ of vertices also labeled $\underline 3$. The
remaining vertices of $\Gamma(S)$, which were left not underlined,
form constant-type rectili\-near paths. Starting from $\xi_3$ and
proceeding clockwise around the center of symmetry, these
paths are first vertical 3-paths (four vertices each) having
alternate set of vertices with constant types 2, 3. The subgraph
induced in $\Gamma(S)$ by the vertices in these 3-paths has parallel
delimiting infinite paths whose edges have the non-vertical
directions alternating. Counterclockwise from $\eta_3$,
infinite constant-type rectilinear paths ascending at an angle of
$30^\circ$ from vertices next to those of $\eta_3$
appear that alternate constant types 1, 3. $\Z_3$-symmetry takes
care of the labels of the other vertices in $\Gamma(S)$. \qfd

\begin{cor}
Both $K_2$-{\it QPDS}s in the proof of Theorem 13 have automorphism group $Z_3$,
yielding asymmetric {\rm(}Penrose{\rm)} tilings of $\Pi$. Its second {\rm(}resp.
first{\rm)} $K_2$-{\rm QPDS} $S$ determines an independent vertex set of $\Gamma(S)$ formed by six {\rm(}resp. an infinite number of{\rm)} hexagons in ${\mathcal P}(S)$ adjacent each in $\Gamma(S)$ to at most {\rm(}resp. more than{\rm)} three hexagons of a common type.
\qfd\end{cor}

The assertion about asymmetric Penrose tilings in Corollary 14 may be compared with a similar assertion in  Section 6 of \cite{Dej} or Section 6 of \cite{DD2}. In the present section, however, as in Sections 2 and 5 above, the triangles of ${\mathcal P}(S)$ may be thought as forming a `mortar' that binds -- and separates -- the hexagons of ${\mathcal P}(S)$ taken as  `bricks', a situation different from that of \cite{Dej,DD2}, even though in both situations  $\delta=3$. This mortar and brick interpretation, however, seems to be restricted by the following conjecture.

\begin{conj}
No $K_2$-{\rm QPDS} in $\Delta$ contains two different triples of pairwise
adjacent hexagons in ${\mathcal P}(S)$ with pairwise different types 1, 2 and 3.
\end{conj}

The following definition starts our approach to classifying the $K_2$-QPDSs in $\Delta$ with no triples of pairwise adjacent hexagons in ${\mathcal P}(S)$ having pairwise different types. Such a classification is completed in Theorem 22, below.

Given a $K_2$-QPDS $S$ in $\Delta$, we say that a {\it hexagonal row} of type $i\in I_3$
for $S$ is a subset ${\mathcal H}$ of hexagons of type $i$ in ${\mathcal P}(S)$,
with ${\mathcal H}$ minimal among the subsets ${\mathcal H}'$ of hexagons of type $i$ in
${\mathcal P}(S)$ having induced automorphism groups
${\mathcal A}({\mathcal H}')=\{f_k:\Delta\rightarrow\Delta ;k\in\Z\}$,
where $f_k$ is the parallel translation given by $f_k(x_1,x_2,x_3)=(x_1+6k,x_2,x_3-6k)$,
for each $(x_1,x_2,x_3)\in\Delta$.
For example, on the left of Figure 6
the hexagons $A_1$ and $A_2$ determine a hexagonal row ${\mathcal H}_A$ of type 2.
Above it, hexagonal rows ${\mathcal H}_B,{\mathcal H}_C,{\mathcal H}_D$ of respective types 1, 1, 3 can be determined by the hexagons $B_1$ and $B_2$, $C_1$ and $C_2$, $D_1$ and $D_2$, respectively.
We may say here that
two `contiguous' hexagonal rows of type 1, (${\mathcal H}_B$ and ${\mathcal H}_C$), are `bordered below' by a hexagonal row of type 2, (${\mathcal H}_A$), and `above'
by a hexagonal row of type 3, (${\mathcal H}_D$). This `sandwiched'
disposition may be continuated, leading to a bisequence of such
hexagonal rows, which can be described by the upward
bisequence $\ldots,2,1,1,3,1,1,2,1,1,3,1,1,2,1,1,3,\ldots$, (for the types $\ldots,2,1,1,3,\ldots$, etc. of the hexagonal rows $\ldots,{\mathcal H}_A,{\mathcal H}_B,{\mathcal H}_C,{\mathcal H}_D,\ldots$, etc., respectively), or
any other bisequencing of substrings $2,1,1,$ and $3,1,1$, (such as $\ldots,2,1,1,2,\ldots$, etc., for the left of Figure 6).
This takes us to the existence of an infinite number of such
{\it sandwiched} $K_2$-QPDS in $\Delta$, as in the following
theorem. We say that the bisequence just written is a
$1^2$-{\it interspersion} of the bisequence
$\ldots,2,3,2,3,2,3,\ldots$, because it is obtained from it by
interspersing double 1's between each two of its terms.

\begin{thm}
For any bisequence $\xi$ whose terms are types $2$ and $3$, there
exists a sandwiched $K_2$-{\rm QPDS} $S$ in $\Delta$ with associated
bisequence obtained by a $1^2$-interspersion of $\xi$.
\end{thm}

\proof This follows from the discussion previous to the statement.
Observe that Theorem 11 is a particular case of this result. \qfd
\begin{figure}
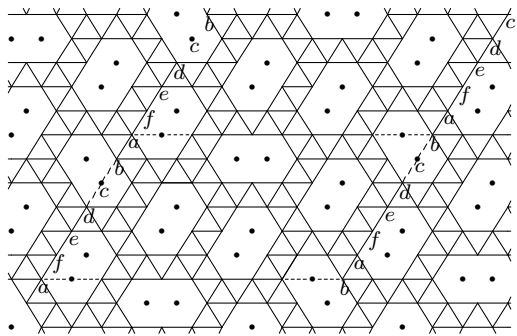

\unitlength=0.40mm
\special{em:linewidth 0.4pt}
\linethickness{0.4pt}

\caption{Illustration for Corollary 17}
\end{figure}
\begin{cor}
If a bisequence $\xi$ as in Theorem 16 has period
$23$, then there exists a graph $\Delta_{m,n}$ containing a
$K_2$-QPDS $S'$ in $\Delta$ which is a quotient of $S$ via the
corresponding projection map
$\rho_{m,n}:\Delta\rightarrow\Delta_{m,n}$ if and only if $6|m$ and
$50|n$.
\end{cor}

\proof Observe Figure 8, which is tilted an angle of $60^\circ$
with respect to the $K_2$-QPDS represented on the left
of Figure 6, so it should be tilted back for a representation of a
$K_2$-QPDS as in the statement, cor\-res\-ponding to the permutation
(3,1,2) of the hexagon types. In Figure 8, a rhomboidal cutout with
a base 50-edges long horizontally and 6-edges upward in
anti-diagonal can be obtained by continuating the sub-rhomboid with
base (or top) $ab$ by means of horizontal translations to the right. The
successive continuating sub-rhomboids here can be realized by
horizontal translation of respective sub-rhomboids (not traced, for clarity)
with bases (or tops) $bc$, $cd$, $de$, $ef$ and $fa$. \qfd

There are $K_2$-QPDSs $S$ in $\Delta$ with hexagons in ${\mathcal
P}(S)$ only of fixed types $i,j\in I_3$, where $i\neq j$, but those of type $i$
given only in pairs of hexagons adjacent in $\Gamma(S)$ and each such pair
surrounded in $\Gamma(S)$ by hexagons of type $j$. Such a situation is illustrated
in Figure 9 with $i=2$ and $j=1$. A table of types as in the
corresponding $\Gamma(S)$ looks locally as follows:

\begin{figure}
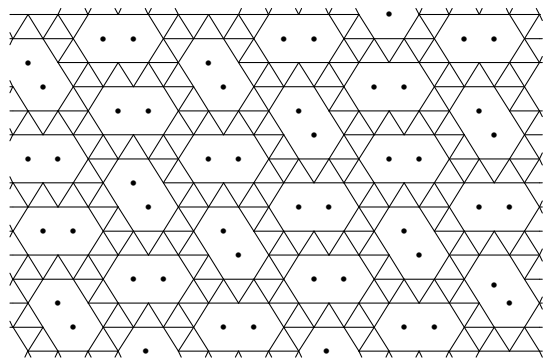

\unitlength=0.40mm
\special{em:linewidth 0.4pt}
\linethickness{0.4pt}

\caption{Pairs of adjacent type 2 hexagons surrounded by type 1 hexagons}
\end{figure}

$$\begin{array}{cccccccccccccc}
^2 &_2 &^1 &_2 &^1 &_1 &^1 &_1 &^2 &_1 &^2 &_2 \\
^1 &_1 &^1 &_1 &^2 &_1 &^2 &_2 &^1 &_2 &^1 &_1 \\
^2 &_1 &^2 &_2 &^1 &_2 &^1 &_1 &^1 &_1 &^2 &_1 \\
^1 &_2 &^1 &_1 &^1 &_1 &^2 &_1 &^2 &_2 &^1 &_2 \\
^1 &_1 &^2 &_1 &^2 &_2 &^1 &_2 &^1 &_1 &^1 &_1 \\
^2 &_2 &^1 &_2 &^1 &_1 &^1 &_1 &^2 &_1 &^2 &_2 \\
\end{array}$$
In such a table, a type $i$ is {\it low} if it has immediately above again a type $i$.
Then, the three descending type sequences that start at a low type $i$, following the diagonal, vertical and anti-diagonal directions, are periodic with periods forming an ordered triple that we called the {\it triple period} $\tau(i,j)$.

Starting with the example above, then rotating $\Delta\subset\Pi$ an angle of $\pm 120\circ$, and finally
reflecting $\Delta$ on the line in $\Pi$ containing the edge of $\Delta[S]$ in an
hexagon of type $i$, the following triple periods $\tau(i,j)$ are found:
$$\begin{array}{cc}
\tau(2,1)=((1^32^2),(1^2212),(121^22)), & \tau(1,3)=((313^21),(3^2131),(3^31^2)), \\
\tau(3,2)=((2^2323),(2^33^2),(2^2323)), & \tau(2,3)=((3^2232),(3^32^2),(3^2232)), \\
\tau(1,2)=((2^31^2),(2^2121),(212^21)), & \tau(3,1)=((131^23),(1^2313),(1^33^2)).
\end{array}$$
It can be seen that similar triple periods can be obtained
by rotating $\Delta$ an angle of $\pm 120^\circ$ and by
reflecting $\Delta$ on a coordinate axis.
Also, it can be seen that there is a graph $\Delta_{m,n}$
containing a $K_2$-QPDS $S'$ which is a quotient of $S$ as in Figure
9 via the correspondence $\rho_{m,n}$ in Corollary 17 if and only if
$50|m$ and $5|n$.

A path in $\Gamma(S)$ is {\it diagonal}
if its edges are in a diagonal disposition in the sense of our
figures and tables. Similar definitions for {\it anti-diagonal} and
{\it vertical} paths hold. The following two theorems generalize the
discussion above, where instead of working with any pair of types
$i,j\in I_3$, we restrict to $i=2$ and $j=1$, since the remaining
five cases given by rotating $\Delta$ angles of $\pm 120^\circ$
and by reflecting it on a coordinate axis differ from this
case just by elementary symmetries.

\begin{thm} There exists exactly one $K_2$-{\rm QPDS} $S$ in $\Delta$ with hexagons
in ${\mathcal P}(S)$ only of types $1,2$, its type $2$ hexagons
given only in maximal diagonal $t$-paths of $\Gamma(S)$ and the
resulting $(t+1)$-tuples of type $2$ hexagons surrounded by type $1$
hexagons, where $t>0$. Moreover, the hexagons of type $1$ in
${\mathcal P}(S)$ are adjacent to either three or four hexagons of
type $2$ in ${\mathcal P}(S)$.
\end{thm}

\begin{cor} There is a graph $\Delta_{m,n}$ containing a $K_2$-{\rm QPDS}
$S'$ in $\Delta$ which is a quotient of the $K_2$-{\rm QPDS} $S$ of
Theorem 18 via $\rho_{m,n}$ if and only if $(10+20t)|m$ and $5|n$.
\end{cor}

\proof The proofs of Theorem 18 and Corollary 19 follow the
discussion previous to their statement and the ideas in the proofs
of Theorem 16 and Corollary 17. Observe that Corollary 17 yields
$K_2$-QPDSs in $\Delta_{m,n}$ which are quotients of the $K_2$-QPDS
in Theorem 16 via $\rho_{m,n}$, for $30|m$ and $5|n$, where
$30=10+20t$ with $t=1$. However, because of the further symmetry in
Theorem 16 and Corollary 17 with respect to that of Theorem 18 and
the present corollary, the equivalent condition in Corollary 17
realizes a smaller minimal value of $m$, namely $m=6$. \qfd

If, for the case $t=1$ in Theorem 18 and Corollary 19, we express
its triple triple as $\tau^t=\tau^1=((1^32^2),(1(12)^2),(121^22))$, then we
can generalize for $t=2$ with the triple period
$\tau^2=((1^42^3),(1(12)^3),((12)21^2)),$ and so on for any value of
$t$:
$$\tau^t=((1^{t+1}2^t),(1(12)^t),(2^{-1}(21)^t)).$$
\begin{figure}
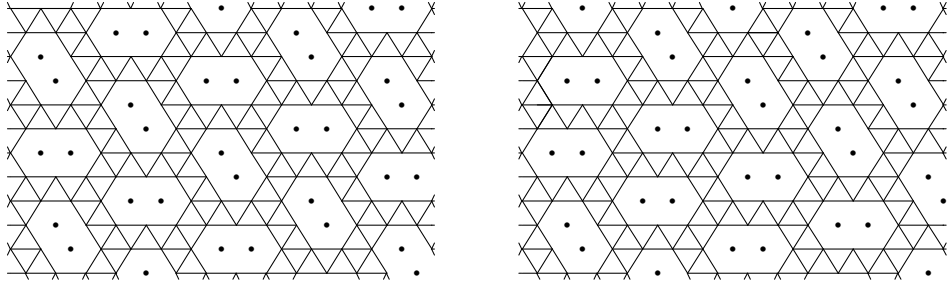

\unitlength=0.40mm
\special{em:linewidth 0.4pt}
\linethickness{0.4pt}

\caption{QPDSs with types 1 and 2 hexagons}
\end{figure}
Now, consider that instead of a $t$-path of type 2 hexagons, we deal
with {\it doubly infinite} paths of type 2 hexagons in ${\mathcal P}(S)$
surrounded by hexagons of type 1. One such a situation is shown
on the left of Figure 10, in which diagonal paths $P_3(2)$ of
${\mathcal P}(S)$ are interspersed with diagonal paths $P_3(1)$,
(where the subindex 3 indicates $P_3(i)$ is a {\it diagonal path} in
$\Gamma(S)$, i.e. with segments between barycenters of hexagons adjacent in $P_3(i)$
perpendicular to {\it anti-diagonals} in $\Delta$, thus the subindex
3, associated with anti-diagonals; and where the type 1 or 2 of the
composing hexagons is set between parentheses). Such a
$K_2$-QPDS will be denoted $S_3(1,2)$, expressing the periodicity of
contiguity and alternation of paths $P_3(2)$ and $P_3(1)$. The
triple period $\tau$ of such an $S_3(1,2)$
could be generalized to the form $S_3(\eta)$, where
$\eta$ is any bisequence formed by types 1 and 2, as on the left of Figure 10, in which part of the complement of the $K_2$-QPDS
$S_3(1^22^2)$ is shown. Observe that periodicity here is not
forced, as nonperiodic bisequences are realizable. These bisequences
can be associated with doubly infinite $\{0,1\}$-sequences, as was
done in Theorem 1 and Corollary 2 of \cite{Dej}. (For the case of
periodic $\eta$, we have Corollary 21 for $K_2$-QPDSs on toroidal
graphs,  below).

\begin{thm} The family of $K_2$-{\rm QPDS}s $S_3(\eta)$ in ${\Delta}$, where $\eta$
varies in the set of doubly infinite $\{1,2\}$-sequences, is in
one-to-one correspondence with the set of points in the real
interval $[0,1]$.
\end{thm}

\proof We modify the representations of QPDSs in $\Delta$ by projecting $\Pi$ onto $\RR^2$, so we can take the $\Delta$-coordinates $x_1,x_2$ in Figure 1 as orthogonal. In this perspective, each hexagon of ${\mathcal P}(S)$ will be displayed with the vertices of its delimiting 6-hole labeled with its type, so each pair of contiguous hexagons with common type 1 or 2 looks like:
$$\begin{array}{lll}
^{111}_{1111}   & & ^{22}_{222} \\
\,\,^{111} & &
\,\,^{222}_{\,\,\,22} \\
\end{array}$$
In this disposition, vertical sequencing represents vertex adjacency, in the representation of Figure 1, along straight paths
in $\Delta$ at an angle of $\pi/6$ over the horizontal paths. For
example, the $K_2$-QPDS $S_3(\eta)$, where $\eta$ has period $(12)$
(of length $I=2)$, becomes representable as in the following
display:

$$\begin{array}{l}
^{\ldots 22222111122222111111\ldots}_{\ldots 122222111111222221111\ldots} \\
^{\ldots 111122222111122222111111\ldots}_{\ldots 2111122222111111222221111\ldots} \\
^{\ldots 2211111122222111122222111111\ldots}_{\ldots 22222111122222111111222221111\ldots}
\end{array}$$
If the horizontal and vertical (or, in $\Delta$, diagonal) periods
created in such a display are indicated by $\phi(I)$ and $\psi(I)$,
respectively, where $I$ is the length of the period of $\eta$, then
the display shows that $(\phi(I),\psi(I))=(\phi(2),\psi(2))=(20,5)$.
A similar representation corresponds for any $S_3(\eta)$ as in the
statement. \qfd

\begin{cor} There exists a graph $\Delta_{m,n}$ containing a $K_2$-{\rm QPDS}
$S'_3(\eta)$ that behaves as a quotient of a $K_2$-{\rm QPDS} $S_3(\eta)$
in $\Delta$ via $\rho_{m,n}$, where $\eta$ is a bisequence as in
Theorem 20 with period $1^{i_1}2^{j_1}1^{i_2}2^{j_2}\ldots
1^{i_t}2^{j_t}$, if and only if $m|\phi(I)$ and $n|\psi(I)$, with
$I=i_1+\ldots+i_t+j_1+\ldots+j_t$ and $(\phi,\psi)(2)=(20,5)$,
$(\phi,\psi)(I)=(10I,5I)$ if $I$ is even $>2$ and
$(\phi,\psi)(I)=(10I,10I)$ if I is odd $>2$.
\end{cor}

\proof In the final display of the proof of Theorem 20, each
horizontal line shows a period composed by the types 1 and 2, and
the last shown line repeats the first one, making it clear that
$(\phi,\psi)(2)=(20,5)$. The remaining cases in the statement are
managed similarly. \qfd

We gather all the results on $K_2$-QPDSs in $\Delta$ obtained above in the following theorem. A $K_2$-QPDS in $\Delta$ is {\it parallel} if it is
as in Theorem 8. Say that a $K_2$-QPDS in $\Delta$ is $\xi${\it
-sandwiched} if it is as in Theorem 16, for some bisequence $\xi$ on
two types in $I_3$; $t${\it -linear} if it is as in Theorem 18, up
to symmetry, for some $t>0$; an $S_i(\eta)$ if it is as in Theorem
20, where $i\in I_3$ and $\eta$ is a bisequence on the types of
$I_3\setminus\{i\}$. We remark that an $(i,j)$-surrounded $K_2$-QPDS
is $1$-linear.

\begin{thm} Up to symmetry,
a $K_2$-{\rm QPDS} in $\Delta$ not as in Theorem 13 or Corollary 15 is
either {\bf (a)} parallel, or {\bf(b)} $t$-linear {\rm(}$t>0${\rm)}, or
{\bf(c)} $\xi$-sandwiched, or {\bf(d)} an $S_i(\xi)$, where $\xi$ is a
bisequence on the two types $\neq i$. \qfd
\end{thm}

It remains to identify completely for which values of integers $m,n$
there exist graphs $\Delta_{m,n}$ containing $K_2$-QPDSs which are
quotients of the $K_2$-QPDSs classified in Theorem 22, though some
of those integer pairs were determined as corollaries above.

Other cases of interest to consider in $\Delta$ and its toroidal quotients are $H$-QPDSs where $H$ is a finite path of length larger than 1.

\end{document}